\documentclass{article}

\usepackage{geometry}

\usepackage{type1cm}        
\usepackage{makeidx}         
\usepackage{graphicx}        
\usepackage{multicol}        
\usepackage{multirow}
\usepackage[bottom]{footmisc}

\usepackage{newtxtext}       %
\usepackage{newtxmath}       

\newcommand{\dd}{\mathrm{d}}
\usepackage{todonotes}

\usepackage{amsmath,amssymb}
\usepackage{nicefrac}
\usepackage{amsfonts}
\newtheorem{definition}{Definition}

\newtheorem{proposition}{Proposition}

\makeindex             

\begin{document}

\title{From kinetic to macroscopic models and back}

\author{
	M. Herty 
	\thanks{
		Institute f\"{u}r Geometrie und Praktische Mathematik -
		RWTH Aachen University --
		Templergraben 55, 52062 Aachen, Germany --
		{\sl herty@igpm.rwth-aachen.de}
	}
	\and
	G. Puppo
	\thanks{
		Dipartimento di Matematica -
		``La Sapienza'' Universit\`{a} di Roma --
		Piazza Aldo Moro 5, 00185 Roma, Italy --
		{\sl gabriella.puppo@uniroma1.it}
	}
	\and
	G. Visconti
	\thanks{
		Institute f\"{u}r Geometrie und Praktische Mathematik -
		RWTH Aachen University --
		Templergraben 55, 52062 Aachen, Germany --
		{\sl visconti@igpm.rwth-aachen.de}
	}
}

\maketitle

\abstract{We study kinetic models for traffic flow characterized by the property of producing backward propagating waves. These waves may  be identified with the phenomenon of stop-and-go waves typically observed on highways. In particular, a refined modeling of the space of the microscopic speeds and of the interaction rate in the kinetic model allows  to obtain weakly unstable backward propagating waves in dense traffic, without relying on  non-local terms or multi--valued fundamental diagrams. A stability analysis of these waves is carried out using the Chapman-Enskog expansion. This leads to a BGK-type model derived as the mesoscopic limit of a Follow-The-Leader or Bando model, and its macroscopic limit belongs to the class of second-order Aw-Rascle and Zhang models.}

\section{Introduction} \label{sec:intro}

There are mainly three modeling scales in the mathematical description of vehicular traffic flow. The microscopic scale is based on the prediction of  trajectories of individual vehicles by systems of ordinary differential equations. The macroscopic scale is based on the assumption that traffic flow behaves like a fluid where individual vehicles cannot be identified, but a macroscopic conservation law for the number of vehicles rules the dynamics. Here, the flow is represented by a density function and evolves in space and time by transport equations. The intermediate scale is the mesoscopic scale. Here, kinetic equations govern the dynamics. Those equations are characterized by a statistical description of the microscopic states of vehicles but, at the same time,  still provide the macroscopic aggregate representation of traffic flow, linking collective dynamics to  interactions among vehicles at a smaller microscopic scale.

In the present chapter we study non-homogeneous kinetic models for vehicular traffic flow. In particular, we investigate the common and well-established idea that non-local terms are necessary to observe backward propagation of waves in dense traffic~\cite{klar1997Enskog}. We show that the model in \cite{PSTV2} naturally encloses backward propagating waves, although these waves may be unstable. We include a first stabilization term  including the effect of uncertainty in the braking rate~\cite{TesiRonc}. We propose a more refined choice of the interaction rate which allows us to obtain weakly unstable waves propagating back in congested traffic situations without considering non-local  terms. More precisely, drawing inspiration from the Knudsen number in kinetic gas-dynamics, we prescribe the interaction rate as a suitable function of the density and its space-derivative. The backward propagating waves may still be unstable in the sense that they may exhibit unbounded growth in time. We study the appearance of these instabilities by considering BGK-type (Bhatnagar, Gross and Krook~\cite{BGK1954}) models in the limit of constant but sufficiently small interaction rates. In this regime   it was been shown in~\cite{BorscheKlar2018} that Enskog-like terms provide a stabilization effect. However, in that work the stabilization  is unfortunately too strong and it implies that for example stop-and-go waves will {\bf not} occur. Following the approach introduced in~\cite{HPRV_BGK2020}, we derive a weakly-unstable BGK model modifying the design of the space of microscopic speeds. Further, we obtain by suitable limits from this mesoscopic representation a  microscopic follow-the-leader~\cite{FTL1961} or Bando~\cite{Bando1995} model, and a macroscopic Aw-Rascle \cite{aw2000SIAP} and Zhang~\cite{Zhang2002} type model.

The chapter is organized as follows. In Section~\ref{sec:2} we introduce  Boltzmann-like kinetic models for traffic flow characterized by binary interactions with over-braking, and we provide an experimental evidence of the backward propagation of waves in dense traffic. In Section~\ref{sec:3} we analyze the stability of these waves by a Chapman-Enskog expansion of the BGK approximation of the full kinetic model, and we compare the results with the Chapman-Enskog expansion of the BGK model in~\cite{BorscheKlar2018} and of the Aw-Rascle and Zhang model. Finally, in Section~\ref{sec:4} we derive a modified version of the BGK model, as in~\cite{HPRV_BGK2020}, and analyze the stability in the case of interactions with over-braking. In Section~\ref{sec:conslusion} we discuss results and future perspectives.

\section{Backward propagation of waves in a kinetic traffic model} \label{sec:2}

A kinetic traffic model for the mesoscopic scale reads as follows
\begin{equation}\label{eq:generalKineticModel}
\partial_t f(x,v,t) + v \partial_x f(x,v,t) = \frac{1}{\varepsilon} Q[f,f](x,v,t),
\end{equation}
where $f(x,v,t)\, : \mathbb{R} \times [0,V_M] \times \mathbb{R}^+ \rightarrow \mathbb{R}^+$ is the mass distribution function of the flow and the 
local traffic  density $\rho(t,x)$ is given by 
\begin{equation}
\rho(x,t)= \int_0^{V_M} f(x,v,t) \dd v.
\end{equation}
We suppose that the space of possible microscopic speeds of the vehicles is bounded  by zero and  a maximum speed $V_M$. Further, we assume that $f(x,v, t=0)$ is such that  density is  limited by a maximum density $\rho_M=\int f(x,v,0) dv < \infty$. Throughout this work, we consider dimensionless quantities and normalize for simplicity  $V_M=1$ and $\rho_M=1$. The source term in~\eqref{eq:generalKineticModel} is commonly called collision kernel, in analogy to kinetic models for gas-dynamics, and it models the change of $f$ due to the microscopic interactions among vehicles. $Q[f,f]$ can be modelled as a non-linear integral operator, typical of Boltzmann-type kernels, or as a linear operator, typical of BGK-type kernels. The quantity $\varepsilon$ is positive, and yields a relaxation rate weighting the relative strength between the convective term and the source term. It is related to the Knudsen number in fluid dynamics. Generally,  $\varepsilon$ can be a function of  density $\rho$, and possibly of its spatial derivative. Here, we consider both the case $\varepsilon=\varepsilon(\rho, \partial_x\rho)$ and the case of a constant  rate $\varepsilon.$

\subsection{A Boltzmann-type kinetic model for traffic flow}

In the collision operator we model the adaptation of vehicles's speeds by binary car--to--car interaction. This behavior is typical for real--world traffic where usually a driver reacts to the actions of the vehicle in front. To describe the interactions we split the 
  operator $Q[f,f]$  in the difference between a gain term and a loss term. The former accounts for the increase of $f(x,v,t)$ when a  vehicle with velocity $v_*$ interacts with a leading vehicle with speed $v^*$, emerging with speed $v$ as a result of the interaction. The latter accounts for the decrease of $f(x,v,t)$ if a  vehicle with velocity $v$ interacts with a  vehicle with speed $v^*$, emerging with speed different from $v$ as a result of the interaction. We assume that the velocity of the leading vehicle remains always unchanged. More specifically,
\begin{equation} \label{eq:collisionKernel}
\begin{aligned}
Q[f,f](x,v,t) =& \int_0^1 \int_0^1 \mathcal{P}(v_*\to v | v^*;\rho) f(x,v_*,t) f(x,v^*,t) \dd v_* \dd v^* \\ &- f(x,v,t) \int_0^1 f(x,v^*,t) \dd v^*.
\end{aligned}
\end{equation}

The core of a kinetic model is the definition of the operator $\mathcal{P}(v_* \to v | v^*;\rho)$ that prescribes, in a probabilistic way, the resulting speed of a vehicle after  interacting with a leading vehicle. The kinetic model for traffic flow studied here is based on the following interaction rules:
\begin{equation} \label{eq:rule}
\mathcal{P}(v_* \to v | v^*;\rho) =
\begin{cases} 
P(\rho) \, \delta_{\min\{ v_* + \Delta_a, V_M \}}(v) + (1-P(\rho)) \, \delta_{\max\{ v_*- \Delta_b, 0 \}}(v) & v_* \le v^* \\ 
P(\rho) \, \delta_{\min\{ v_*+ \Delta_a, V_M \}} (v) + (1-P(\rho)) \, \delta_{\max\{ v^*- \Delta_b, 0 \}} (v) & v_* > v^* 
\end{cases}
\end{equation}
where $P(\rho) \in [0,1]$ is a decreasing function of the density modeling the  probability of accelerating. The parameters $\Delta_a$ and $\Delta_b$ are the acceleration and the braking parameters, respectively, where  $\Delta_a$ is the instantaneous physical acceleration of a vehicle. The parameter $\Delta_b$ instead corresponds to an uncertainty in the estimate of the other vehicle's speed. Indeed, $\Delta_b=0$ corresponds to no uncertainty: the  vehicle has an exact perception of velocities, and therefore is able to maintain its own speed $v = v_*$ when it interacts with a faster  vehicle (i.e. when $v_*<v^*$), while it can brake exactly to the speed $v = v^*$ in case a slower  vehicle is ahead (i.e. when $v_*>v^*$). For $\Delta_b=0$  the model ~\cite{PSTV2} is recovered. More details on the case $\Delta_b>0$ can be found in \cite{TesiRonc}. Note that the model is continuous across the line $v_*=v^*$, ensuring well-posedness, see \cite{PgSmTaVg3}, and that mass conservation holds:
$$
	\mathcal{P}(v_*\to v | v^*;\rho) \geq 0, \quad \int_0^1 \mathcal{P}(v_*\to v | v^*;\rho) \dd v = 1.
$$

In the space homogeneous case $f = f(v,t)$, the model~\eqref{eq:generalKineticModel} reduces to a relaxation to equilibrium which is characterized by a function $M_f(v;\rho)$ such that $Q[M_{f},M_{f}]=0$. In analogy to kinetic models for rarefied gas dynamics, the function $M_f$ will be called Maxwellian and it allows us to define the flux and the mean speed of vehicles at equilibrium as
\begin{equation} \label{eq:macroQuantitiesEq}
	F_\text{eq}(\rho) = \left(\rho U_\text{eq}(\rho)\right) = \int_0^1 v M_f(v;\rho) \dd v, \quad U_\text{eq}(\rho) = \frac{1}{\rho} \int_0^1 v M_f(v;\rho) \dd v.
\end{equation}

For $\Delta_b=0$ it is proven, cf.~\cite{PSTV2}, that stable equilibria are uniquely defined by the local density. Moreover, the Maxwellian is a known function of $v$, it can be explicitely computed, and depends on $x$ and $t$ only through the local density $\rho(x,t)$. Further, in the space homogeneous case, the density is a scalar parameter fixed at the initial time. However, unstable equilibria may also occur, for which the Maxwellian depends not only on $\rho$, but also on the initial distribution $f(x,v,t=0)$. These equilibria are unstable under perturbation of the  initial datum. The Maxwellian corresponding to the stable equilibria is a  finite weighted sum of Dirac's functions for any initial distribution. If the braking uncertainty $\Delta_b \neq 0$, it has been shown in~\cite{TesiRonc} that the equilibria corresponding to a given density are unique, and all equilibria are stable. 

\begin{figure}
	\centering
	\includegraphics[width=\textwidth]{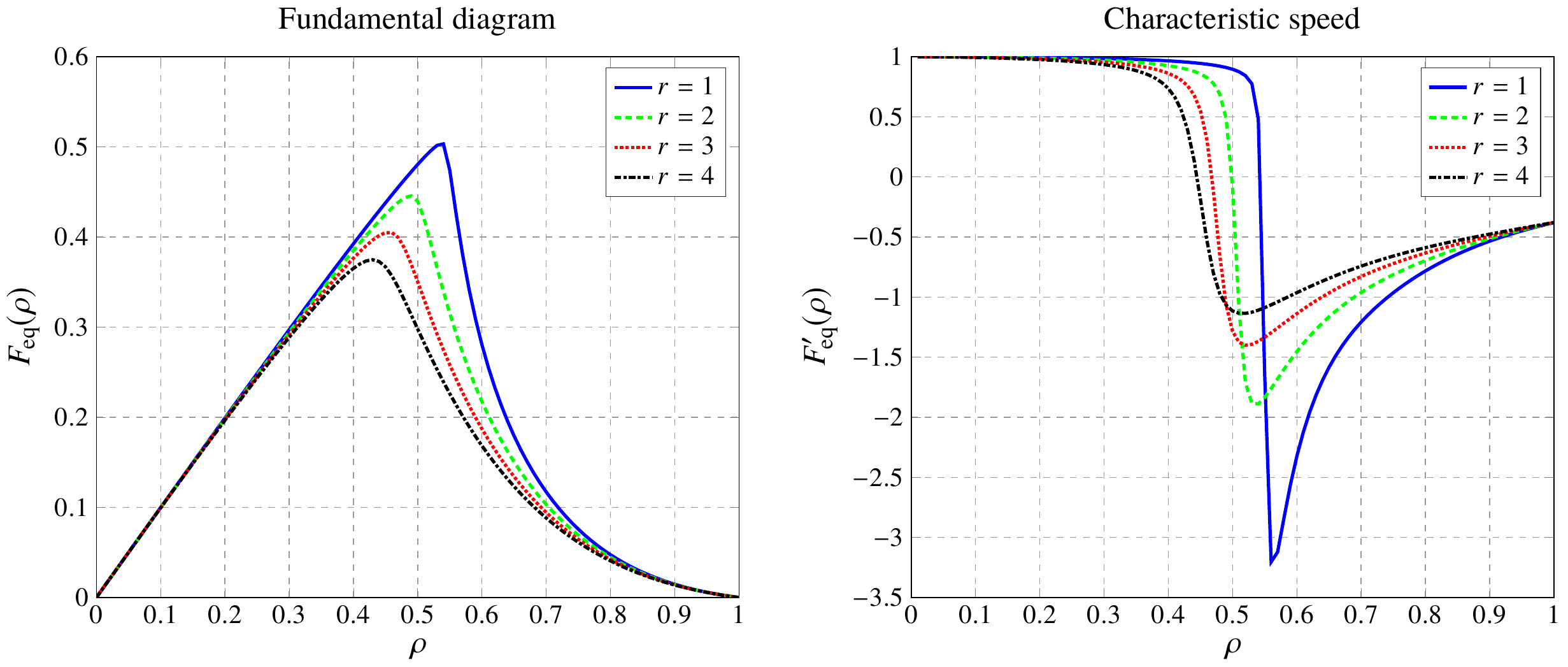}
	\caption{Fundamental diagrams (left) and characteristic speed (right) with $48$ discrete microscopic speeds, $\Delta_a=\frac14$ and $\Delta_b=\frac{\Delta_a}{r}$, $r=1,2,3,4$.}
		\label{fig:fundDiagrams}
\end{figure}

In Figure~\ref{fig:fundDiagrams} we show the equilibrium flux $F_\text{eq}(\rho)$, also known as fundamental diagram, and the characteristic speed $F_\text{eq}^\prime(\rho)$ obtained numerically by using $48$ discrete equidistant discretization points in the velocity phase space, a fixed value of the acceleration parameter $\Delta_a=\frac{V_M}{4}=\frac14$ and different values of the uncertainty $\Delta_b$ such that $r=\frac{\Delta_a}{\Delta_b}=1,2,3,4$. In all cases, the fundamental diagram is characterized by two phases. For low values of $\rho$ the flux is linear in $\rho$. This is the phase of free flow. For larger values of $\rho$, the role of the interactions increases, and the flux decreases. This corresponds to the congested phase of traffic flow. The value of the density for which the change between congested and free flow occurs is called {\em critical density}. Note that the road capacity, i.e. the maximum of the flux, decreases as the uncertainty $\Delta_b$ increases. 

\subsection{Propagation of waves}

Integrating equation~\eqref{eq:generalKineticModel} in velocity space, the right-hand side vanishes because of mass conservation, and one obtains the evolution equation for the density
\begin{equation}\label{eq:macroKin}
\partial_t \rho(x,t) + \partial_x F(x,t;f) = 0, \quad F(x,t;f) = \int_0^1 v f(x,v,t) \dd v,
\end{equation}
where $F$ is the macroscopic flux obtained through the kinetic model. If the system approaches equilibrium, $f \rightarrow M_{f}$, and the macroscopic equation reduces to the equilibrium equation
\begin{equation}\label{eq:macroEq}
\partial_t \rho(x,t) + \partial_x F_\text{eq}(\rho(x,t)) = 0, \quad F_\text{eq}(\rho(x,t)) = \int_0^1 v M_{f}(v;\rho) \dd v.
\end{equation}
Since the Maxwellian is  defined by $\rho$, the equilibrium equation \eqref{eq:macroEq} is closed, and it is a well defined scalar conservation law where the flux function $F_\text{eq}(\rho)$ is the fundamental diagram. On the other hand, when the system is not at equilibrium, the macroscopic equation~\eqref{eq:macroKin} is still coupled to the kinetic equation~\eqref{eq:generalKineticModel}. 

At the mesoscopic scale, the relaxation speed defined by $\varepsilon$ plays a crucial role since, balancing the weight between the convection and the source term, it allows us to define the regimes of the kinetic model. If we allow for $\varepsilon = 0$, i.e. we suppose that the interactions are so frequent to instantaneously relax $f$ to the local equilibrium distribution $M_f$, we are in the so-called equilibrium flow regime where~\eqref{eq:generalKineticModel} reduces to the conservation law for the density~\eqref{eq:macroEq}. Instead, we expect that if $\varepsilon$ is small, but not vanishing, then we are either in a regime where the kinetic equation~\eqref{eq:generalKineticModel} reduces to a {\em perturbed} continuity equation~\eqref{eq:macroEq}  or where  the kinetic equation can be approximated by an extended continuum hydrodynamic system of equations as, for example,  the Aw-Rascle and Zhang model. For $\varepsilon \asymp 1$, but not too large, we are in the kinetic regime  and finally for $\varepsilon \gg 1$ we obtain the regime of the collision-less kinetic equation where the convective term dominates.

In regimes characterized by a small value of $\varepsilon$, we expect that the conservation law \eqref{eq:macroEq} should provide a good approximation to the behavior of the solution; in particular smooth waves should travel along the characteristics given by $\partial_\rho F_\text{eq}(\rho)$. Thus, looking at the right panel of Figure~\ref{fig:fundDiagrams}, we expect that signals move towards the right in the free flow phase and towards the left in the congested flow phase. 

However, in the kinetic regime where $\varepsilon >> 1$ signals should  always propagate towards the right  since the microscopic velocities in traffic are non-negative. This happens also for congested traffic regimes, because the characteristics in the transport term coincide with the microscopic speeds. As observed in~\cite{klar1997Enskog} this can be seen by computing the implicit solution to~\eqref{eq:generalKineticModel}  
$$
f(x,v,t) = f(x-vt,v,t=0) + \int_0^t Q[f,f](x+v(s-t),v,s) \dd s.
$$
The distribution function $f$ at $x$ and $t$ depends only on the distribution function at the values $y \leq x$ and $s \leq t$, since $v$ is non--negative. Thus, apparently, traffic jams in dense flow are not allowed to travel backwards. Several models were introduced in the mathematical literature~\cite{FermoTosin13,klar1997Enskog} trying to overcome this drawback. 

Numerical evidence suggests strongly that this picture is naive, and that the interaction of the source term, given by the collision operator, and the transport term, is more subtle. We observe instead a smooth transition between the solutions of the equilibrium equation, where signals move backward in congested flow, and solutions of the kinetic equation. Here too in fact the propagation speed of smooth waves can be negative.

\begin{figure}
	\includegraphics[width=0.5\textwidth]{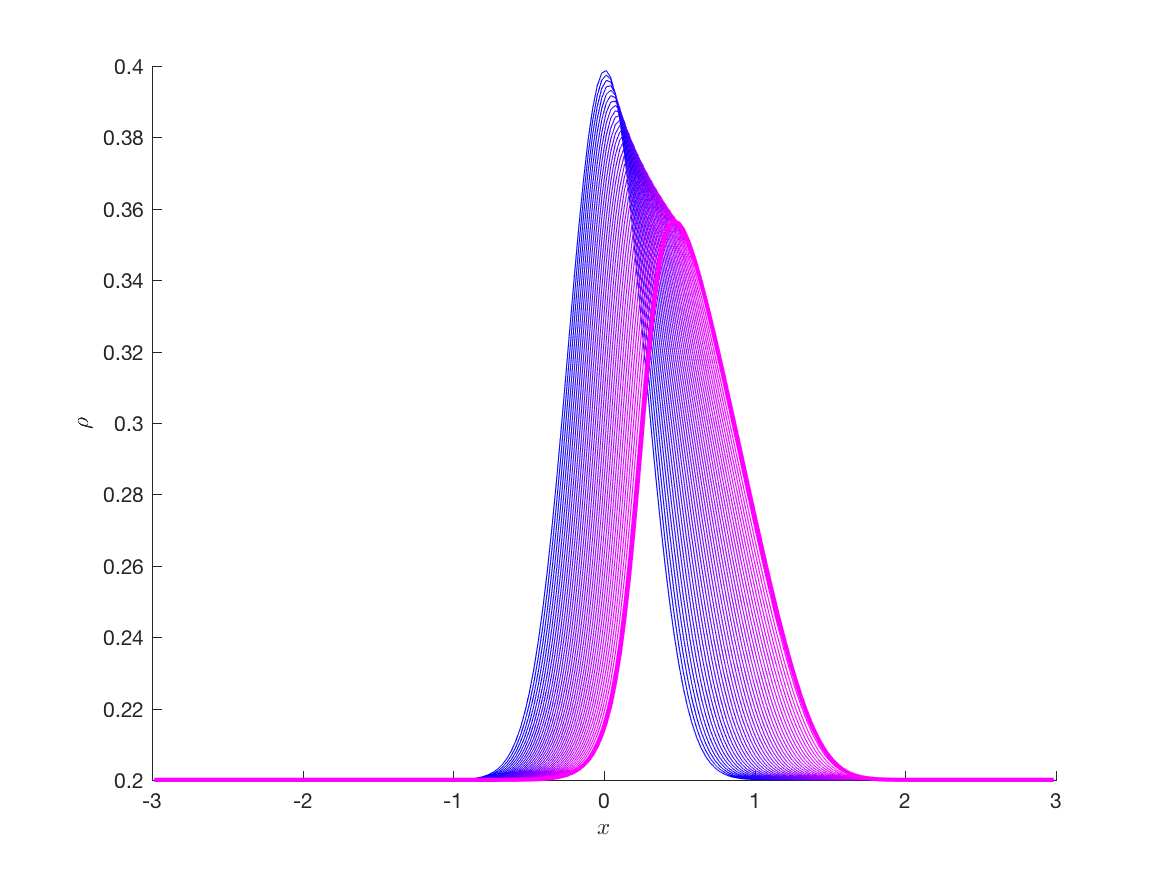}
	\includegraphics[width=0.5\textwidth]{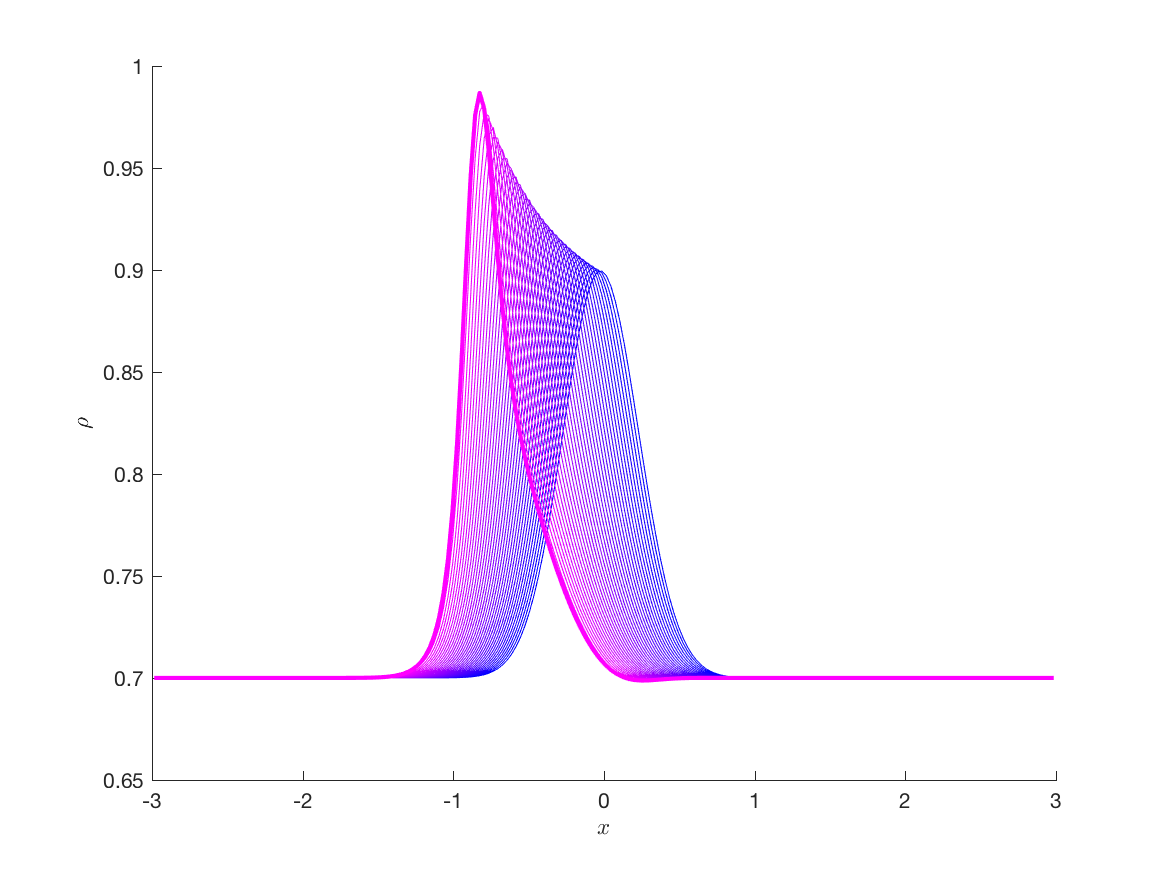}
	\caption{Time evolution of a density bump in the free-flow phase (left) and in the congested-flow phase (right). The initial condition is drawn in blue, and the solution shades towards magenta, as time increases.} \label{fig:kineticEvolution}
\end{figure}

To illustrate this point, we show the evolution of the solution of the kinetic model~\eqref{eq:generalKineticModel} in a few typical cases. In particular, we consider propagating a smooth perturbation in the density
$$
	\rho_0(x) = a + b e^{-8x^2}
$$
and periodic boundary conditions. The initial distribution is Maxwellian. The solution is computed with a first order numerical method, using the local Lax Friedrich's flux. The choice of the numerical flux is crucial: a standard upwind flux, computed following the characteristics of the transport term, would in fact be unstable, in the congested phase, because the direction of the flow {\em does not} coincide with the direction of the characteristics. Since the collision term becomes stiff for small $\varepsilon$, we penalize the collision term with a BGK operator, as in \cite{DimarcoPareschi2013}.

We use $4$ discrete speeds with $\Delta_a=\Delta_b=\frac14$; space is discretized by $200$ cells and the final time is $t_f=1$, while $\varepsilon=0.01$. The solution is shown at different times, starting from the blue curve at $t=0$, and ending with the magenta thick profile, at $t=1$.  In the left panel of Figure~\ref{fig:kineticEvolution}, we take $a=0.2$ and $b=0.2$. The perturbation in the density is below the critical density. Thus, the density profile moves towards the right, as it would occur also in the equilibrium equation. The shape of the initial data is deformed mainly by numerical diffusion, because the flux is almost linear. In the right panel of Figure~\ref{fig:kineticEvolution}, we choose $a=0.7$ and $b=0.2$, so that the initial perturbation has the same amplitude as before, but it occurs on the dense traffic regime. Now, we observe propagation of the wave towards the left, although the characteristics point towards the right. This means that the propagation speed is governed by the interaction between the collision kernel and the transport term, which reproduces the behavior of the fundamental diagram of the equilibrium equation, where indeed we observe negative characteristics. Note that the height of the density peak now {\em increases} with time: the solution has the correct propagation speed, but it is unstable.

\begin{figure}
	\includegraphics[width=0.5\textwidth]{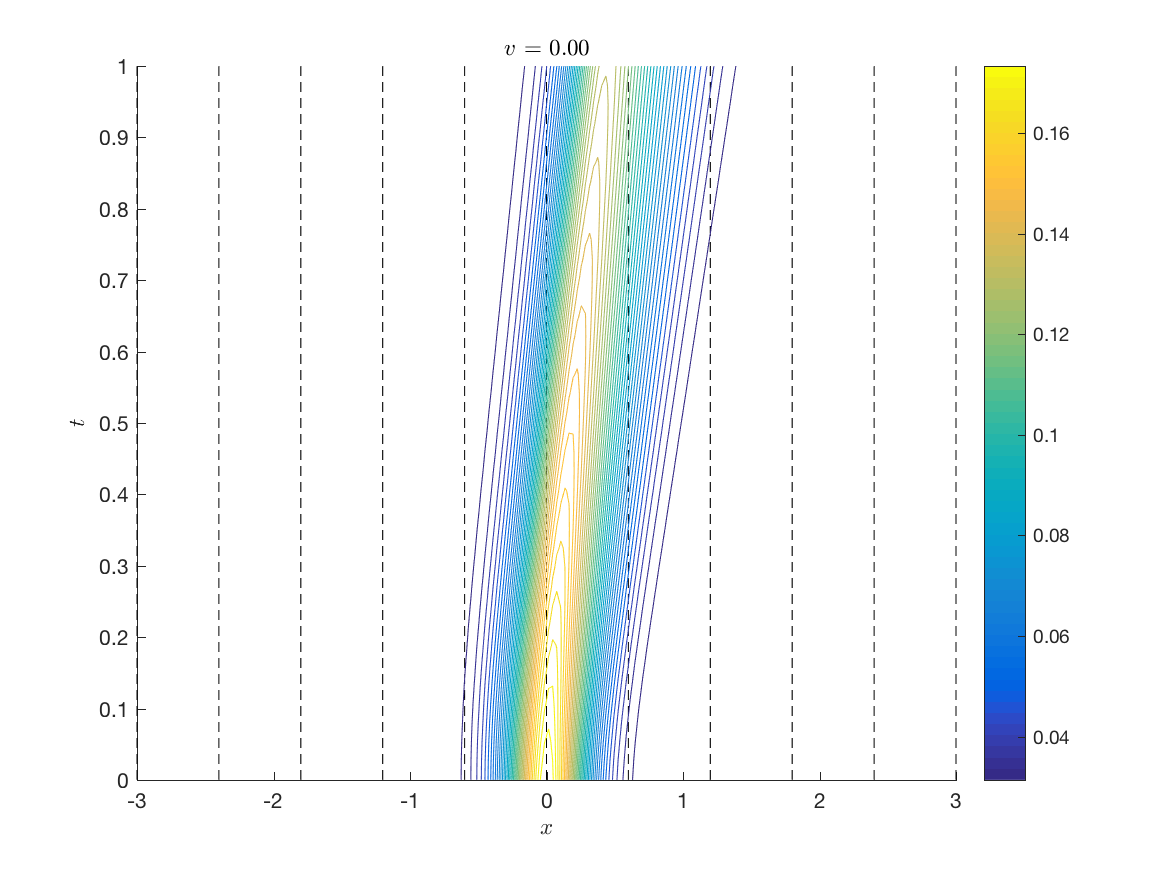}
	\includegraphics[width=0.5\textwidth]{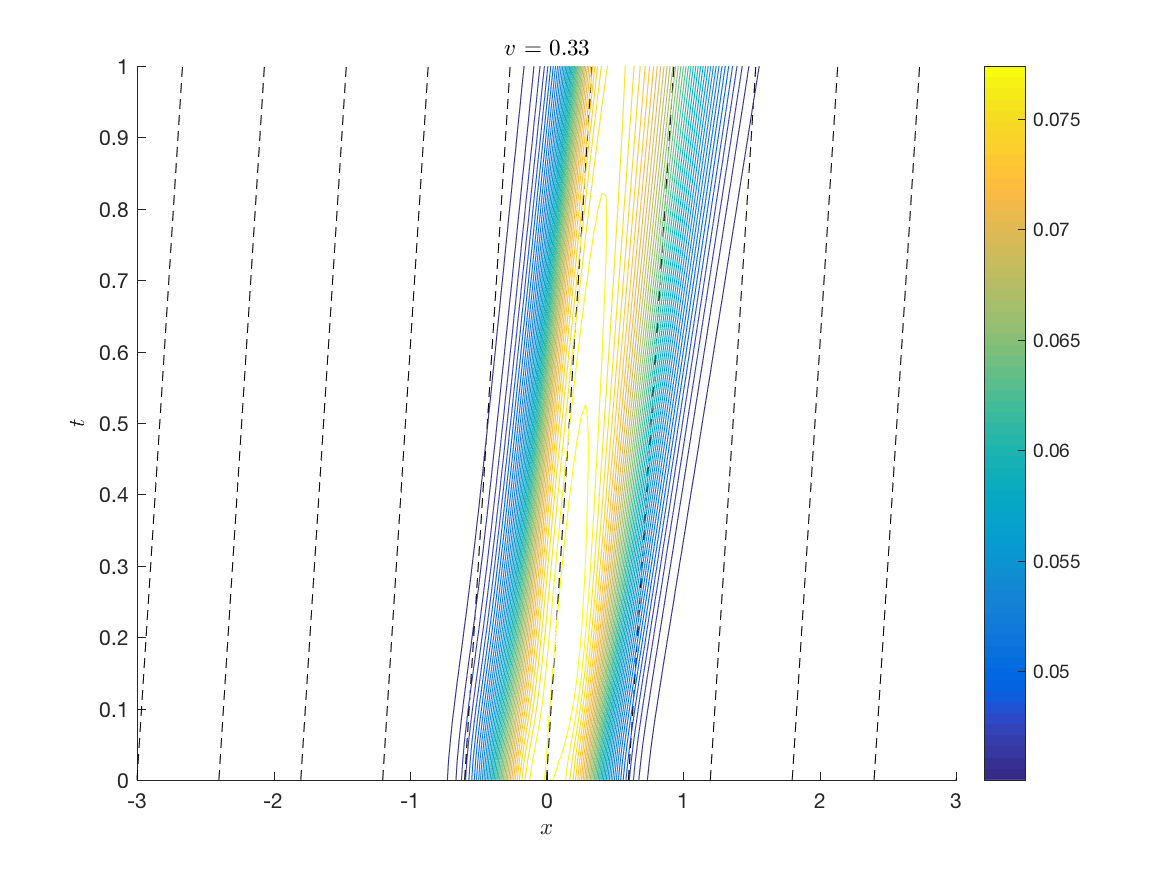}
	\includegraphics[width=0.5\textwidth]{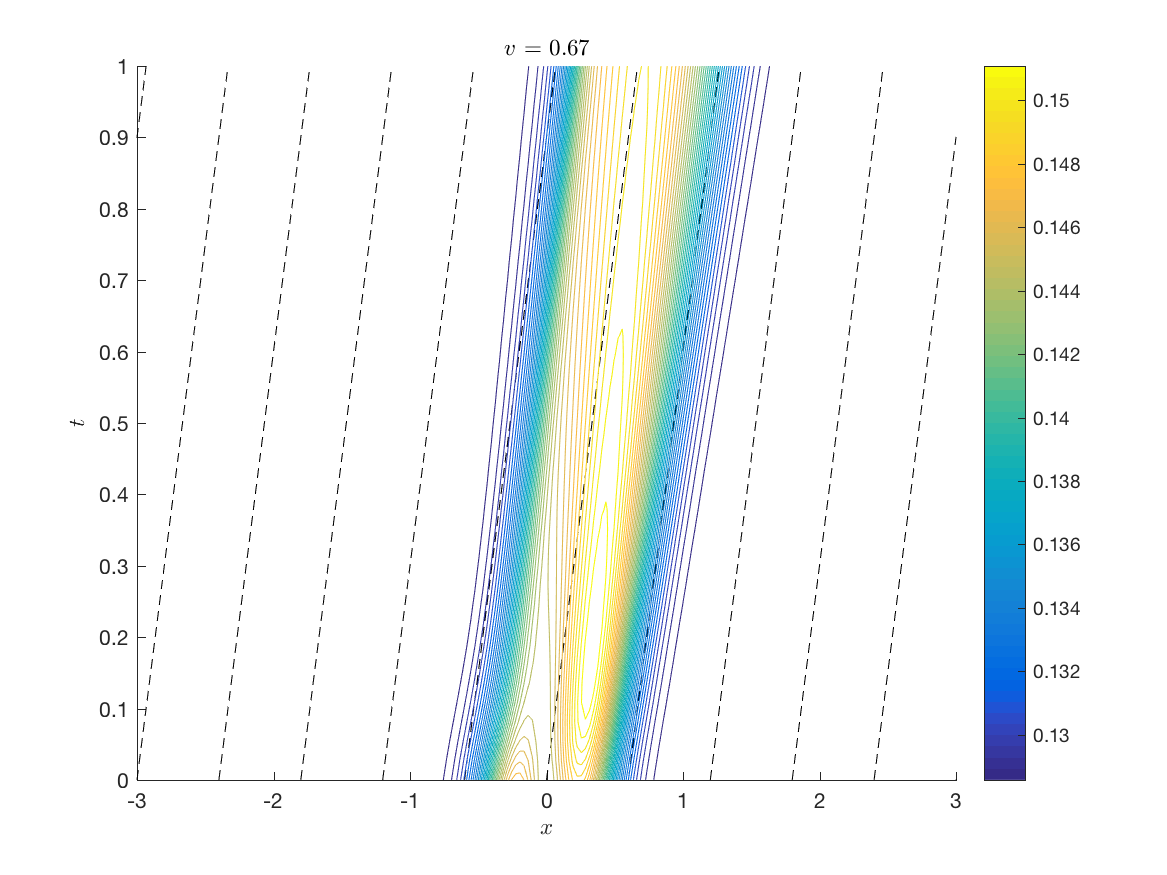}
	\includegraphics[width=0.5\textwidth]{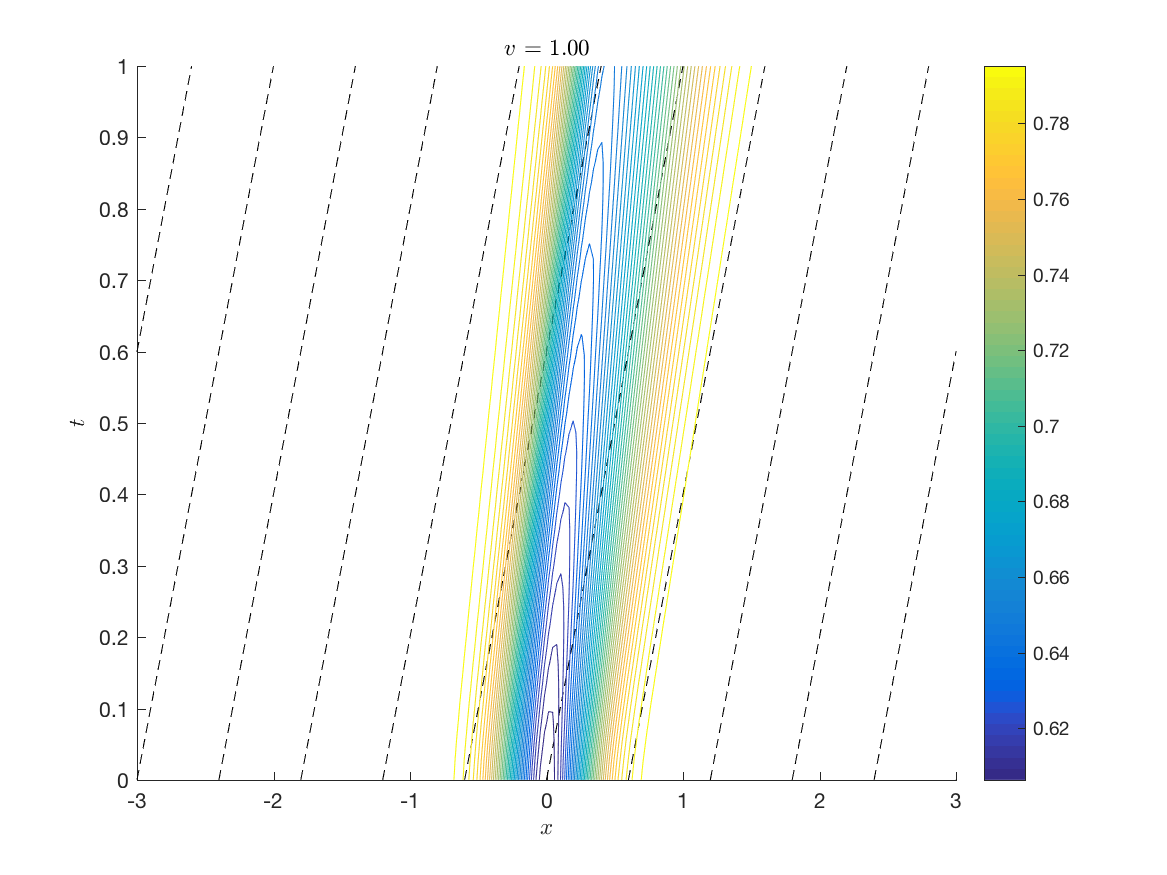}

	\caption{Space-time evolution of the distribution function for each fixed value of the microscopic speed $v$, superposed to the corresponding characteristic speed (black dashed lines), during the time evolution of the density bump in free-flow.\label{fig:charFreeFlow}}
\end{figure}

\begin{figure}[t!]
	\includegraphics[width=0.5\textwidth]{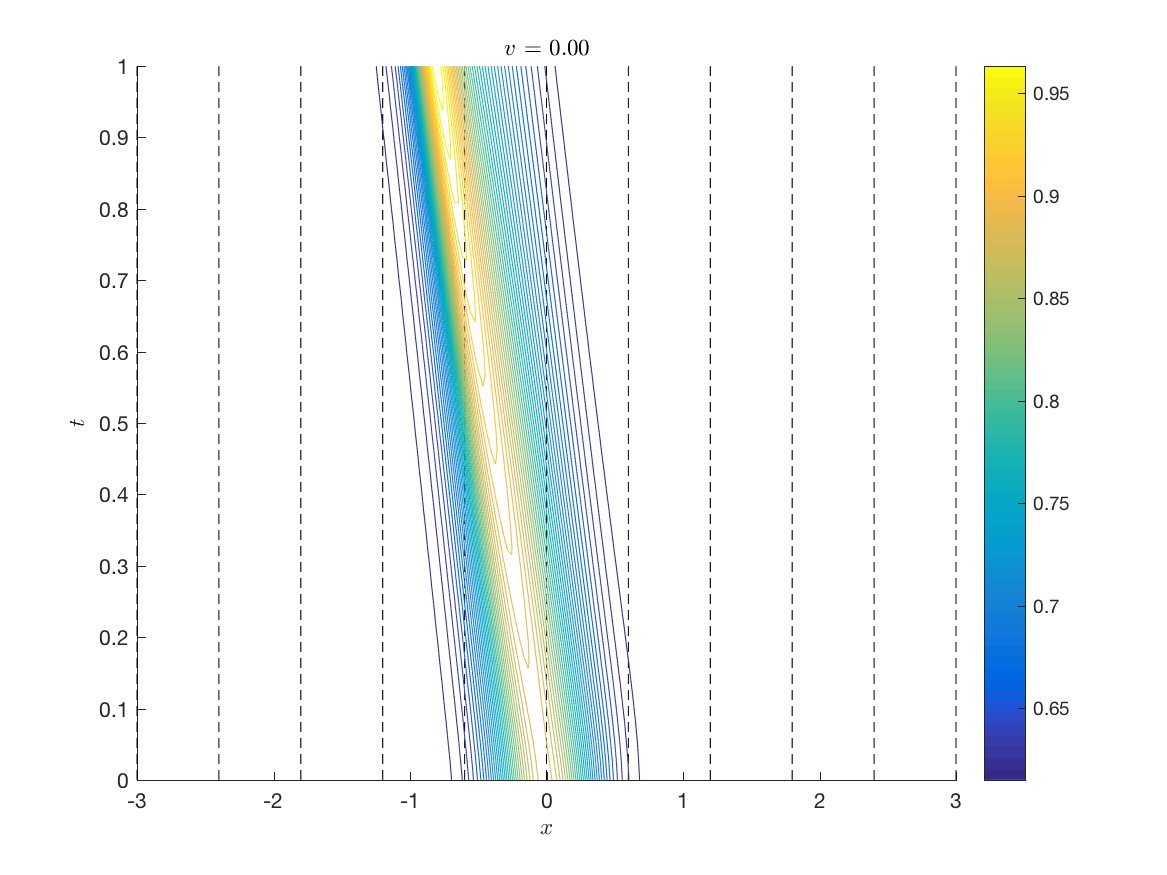}
	\includegraphics[width=0.5\textwidth]{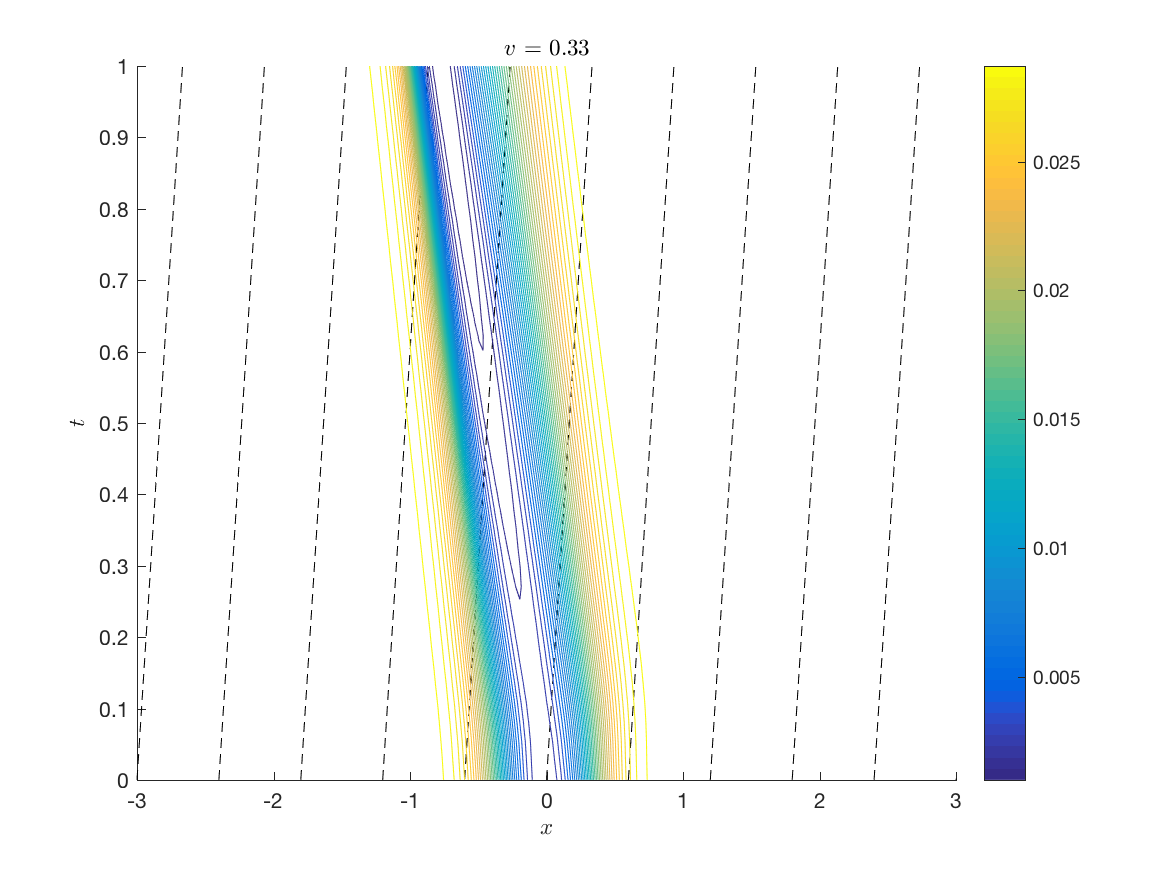}
	\includegraphics[width=0.5\textwidth]{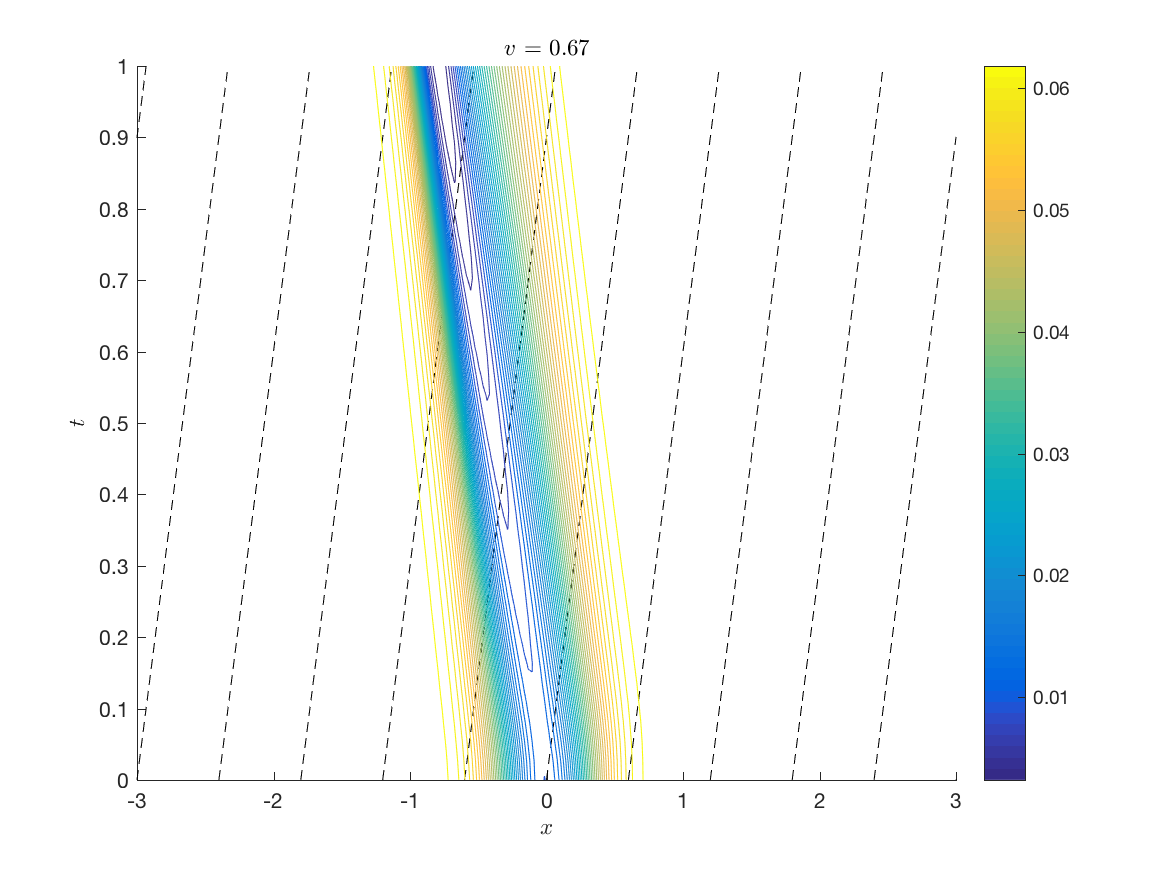}
	\includegraphics[width=0.5\textwidth]{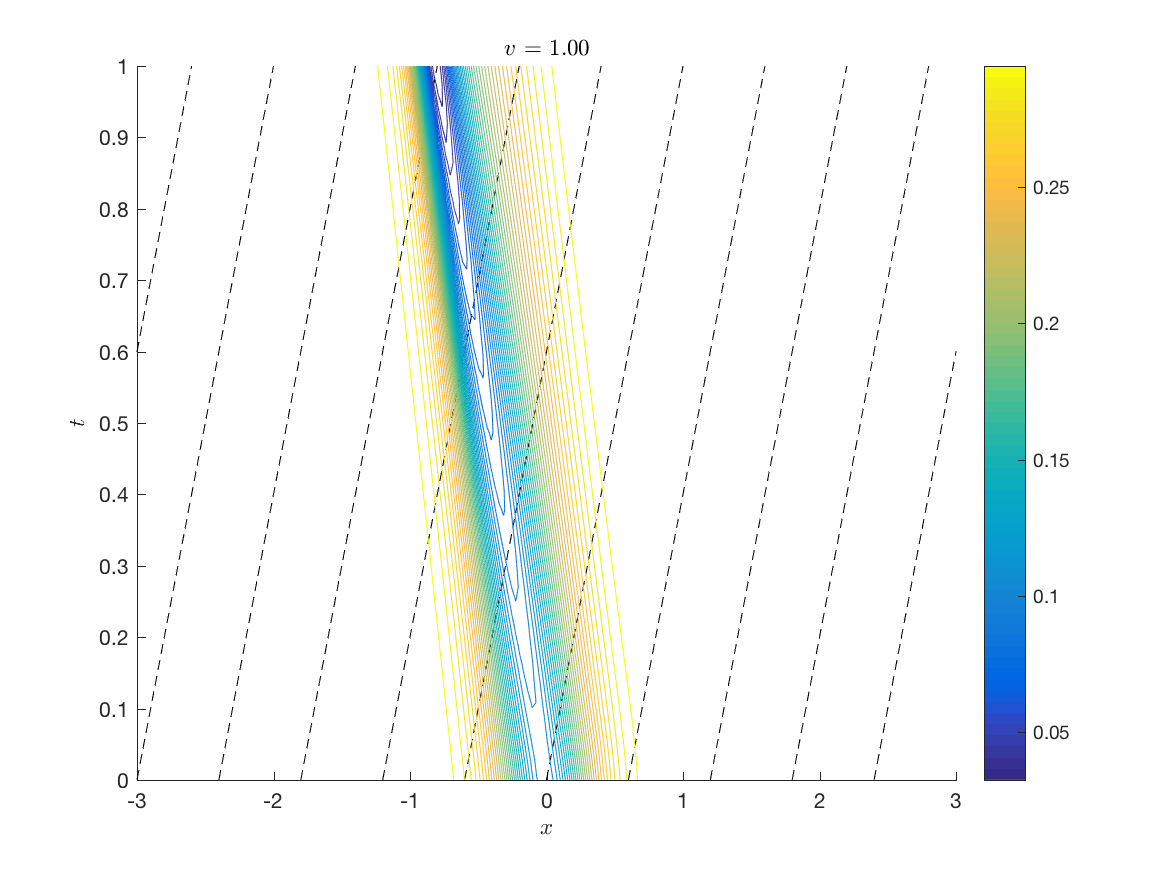}
	\caption{Space-time evolution of the distribution function for each fixed value of the microscopic speed $v$, superposed to the corresponding characteristic speed (black dashed lines), during the time evolution of the density bump in congested-flow.\label{fig:charCongestedFlow}}
\end{figure}

These considerations are further investigated by looking into the distribution function of the kinetic model. We draw contour plots of the space-time behavior of $f$, for each fixed value of the microscopic speed $v$. Since we are considering 4 microscopic velocities, we obtain 4 different plots. In the plots, we compare the time development of the solution $f$ with the corresponding   characteristic speed of the transport term,  drawn with parallel dashed black lines. In the case of the density profile in the free flow phase, we see that the signal propagates towards the right and along characteristics, Figure~\ref{fig:charFreeFlow}. Instead, in the case of the density profile in the congested phase, it is clear that the signal propagates towards the left and {\em across }characteristics. Thus the information on the propagation is contained in the interaction of the collision kernel and the transport term, rather than in the convective term alone.

A constant choice of  $\varepsilon$ however is not satisfactory. In fact, in analogy with the Knudsen number in gas models, $\varepsilon$ should be a decreasing function of the density. In this way, $\varepsilon$ becomes large in the free flow phase since the interactions are less frequent and the convective term rules the dynamics. On the contrary, $\varepsilon$ should become small when the density increases, since the relaxation towards equilibrium should be fast when $\rho$ is high and interactions among vehicles are dominant. Further, we also expect that $\varepsilon$ should decrease  when the traffic thickens, i.e. when $\rho_x$ is large and positive. 
A choice respecting this argument  is 
\begin{equation}\label{eq:VariableEpsilon}
\varepsilon(\rho,\rho_x) = \frac{1}{\max\left\{ \frac{1}{1-\min\left\{\rho,\varepsilon_0\right\}^2},1+\Big(\max\left\{\rho_x,0\right\}\Big)^2\right\}},
\end{equation}
where $\varepsilon_0$ is a threshold to prevent division by zero. The dependence on $\max(\rho_x,0)$ is crucial to prevent overshoots above the maximum density $\rho_M=1$, when the density profile is very steep. This might happen if the  density increases sharply, as when a fast, low density traffic impinges against a slow congested region. In this case, the presence of $\partial_x \rho$ accounts for the need to look ahead. It replaces the non locality of the collision term introduced in~\cite{klar1997Enskog}.

\begin{figure}
	\centering
	\includegraphics[width=0.5\textwidth]{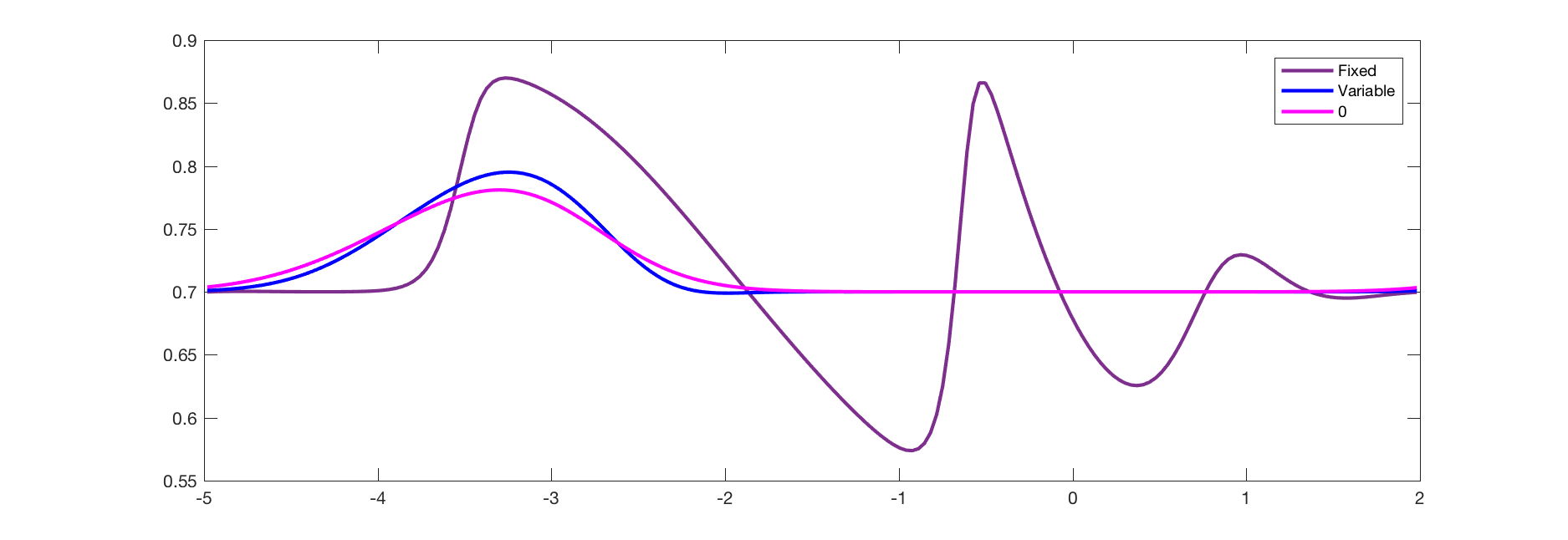}
	\includegraphics[width=0.5\textwidth]{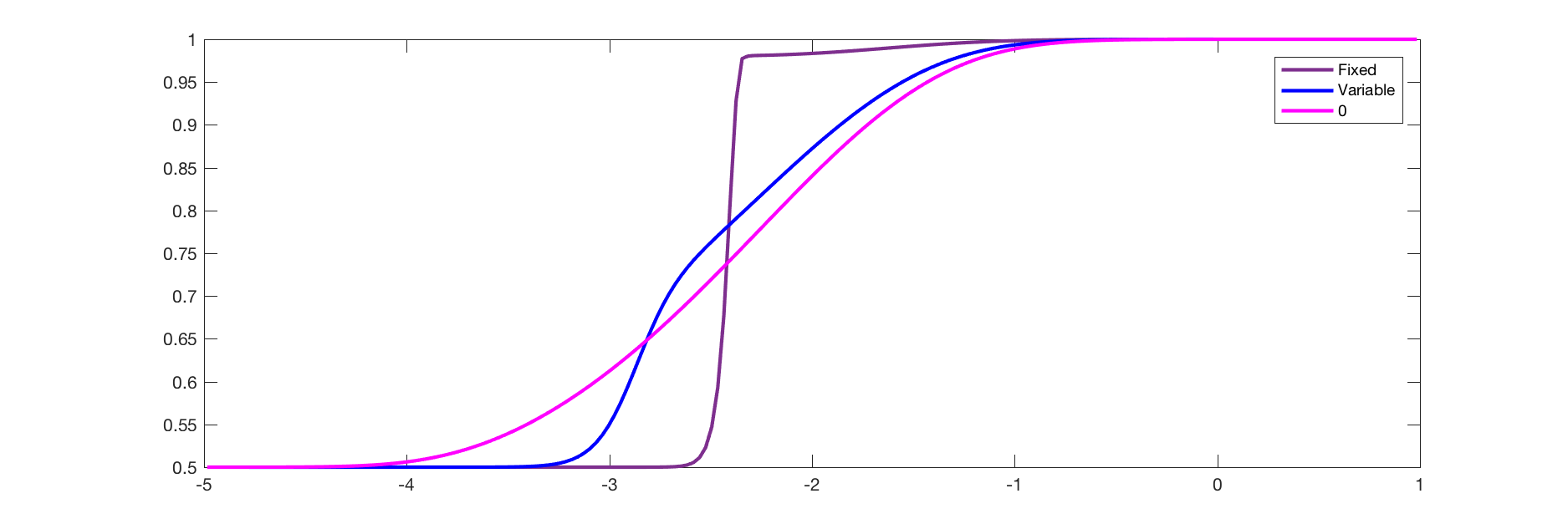}
	\caption{Comparison of kinetic solutions with variable $\varepsilon$, as in eq. eqref{eq:VariableEpsilon}, with fixed $\varepsilon=0.01$, and with the equilibrium solution $\varepsilon=0$. Top: solution with a smooth profile in the congested regime. Bottom: solution of a Riemann Problem, corresponding to a red light.\label{fig:VariableEpsilon}}
\end{figure}

A comparison between a fixed $\varepsilon$ and the variable collision time of \eqref{eq:VariableEpsilon} is shown in Fig. \ref{fig:VariableEpsilon}. The top part of the figure contains the evolution of the high density profile with $a=0.7$ and $b=0.2$ up to time $t=10$. We see that with the variable collision time, the profile propagates to the left, developing waves which resemble stop and go waves. The fixed value of $\varepsilon=0.01$ prevents the developing of these waves, because the relaxation rate is very strong even when the interaction should be weak. As a comparison, we also show the solution obtained with the equilibrium equation \eqref{eq:macroEq}. 

The bottom part of the figure shows the solution obtained for a Riemann problem mimicking a stream of low density traffic impinging against a queue. Here, the kinetic solution with variable $\varepsilon$ develops correctly a shock wave, while the equilibrium solution yields a smooth wave, because, in the congested regime, the fundamental diagram of \eqref{eq:macroEq} is convex.

\section{Analysis of instabilities via Chapman-Enskog expansion} \label{sec:3}


The presence of instabilities is investigated using a formal Chapman-Enskog expansion. For the sake of simplicity, the analysis is performed using the BGK approximation of the Boltzmann-type collision kernel~\eqref{eq:collisionKernel}, for constant but small values of $\varepsilon$. Unstable waves are present also in this linearized setting~\cite{HPRV_BGK2020}. 


\subsection{BGK approximations with and without non-local terms}

The BGK approximation to the kinetic model ~\eqref{eq:generalKineticModel} reads
\begin{equation} \label{eq:BGK}
	\partial_t f(x,v,t) + v \partial_x f(x,v,t) = \frac{1}{\varepsilon} \left( M_f(v;\rho)-f(x,v,t) \right).
\end{equation}
The BGK model is an approximation of the full kinetic equation, which holds for small values of $\varepsilon$. In fact, \eqref{eq:generalKineticModel} and~\eqref{eq:BGK} have, by construction, the same equilibrium solution. This further motivates the use of the BGK approximation to investigate the appearance of instabilities in dense traffic, i.e. in the regime of large densities and small $\varepsilon$.

The Chapman-Enskog expansion allows us to study the behavior of~\eqref{eq:BGK} when $f$ is a first order perturbation in $\varepsilon$ around the equilibrium distribution $M_f(v;\rho)$. In particular, we consider fixed and small values of $\varepsilon$. Then, plugging the expansion
$$
f(x,v,t) =  M_f(v;\rho) + \varepsilon f_1(x,v,t), \quad \text{with } \int_0^1 f_1(x,v,t) \dd v = 0,
$$
into~\eqref{eq:BGK} and integrating with respect to the velocity leads to the advection-diffusion equation
\begin{equation} \label{eq:advdiff}
\partial_t \rho(x,t) + \partial_x F_{\text{eq}}(\rho(x,t)) = \varepsilon \partial_x (\mu(\rho) \rho_x(x,t)),
\end{equation}
where the diffusion coefficient $\mu(\rho)$ is given by 
\begin{equation} \label{eq:diffBGK}
\begin{aligned}
\mu_{\text{BGK}}(\rho) &= \int_0^1 v^2 \partial_\rho M_f(v;\rho) \dd v - \left(\int_0^1 v \partial_\rho M_f(v;\rho) \dd v\right)^2 \\ &= \int_0^1 v^2 \partial_\rho M_f(v;\rho) \dd v - F_\text{eq}^\prime(\rho)^2.
\end{aligned}
\end{equation}

If $\mu(\rho)<0$ then the advection-diffusion equation is ill-posed and therefore may exhibit solutions with unbounded growth.  In the case of the kinetic model~\eqref{eq:BGK}, the sign of the diffusion coefficient depends on the equilibrium distribution $M_f$. The request $\mu_\text{BGK}(\rho) > 0$ is
\begin{equation} \label{eq:stabilityBGK}
	\partial_\rho \left( \int_0^1 v^2 M_f(v;\rho) \dd v \right) > F_\text{eq}^\prime(\rho)^2
\end{equation}
and, since $F_\text{eq}(\rho)$ is the fundamental diagram at equilibrium, this condition requires that the square of the characteristic velocities is bounded by the variation of the kinetic energy in each regime.

Below we recall the result in~\cite{HPRV_BGK2020}, which proves that the instability of the solution does not depend on the choice of the equilibrium distribution and in fact occurs for any suitable equilibrium of kinetic traffic models.

\begin{proposition} \label{th:negativeDiff}
	Assume that $\exists \, \widetilde{\rho} \in (0,1)$ such that
	\begin{equation} \label{eq:conditions}
	F_\text{eq}^\prime(\rho) = \int_0^1 v \partial_\rho M_f(v;\rho) \dd v < 0, \quad \partial_\rho \mathrm{Var}(v) = \partial_\rho \int_0^1 (v-U_\text{eq}(\rho))^2 M_f(v;\rho) \dd v < 0
	\end{equation}
	for all $\rho\in(\widetilde{\rho},1)$. Then the quantity $\mu(\rho)$ given in~\eqref{eq:diffBGK} is negative $\forall\,\rho\in(\widetilde{\rho},1)$.
\end{proposition}

\begin{figure}[t!]
	\centering
	\includegraphics[width=\textwidth]{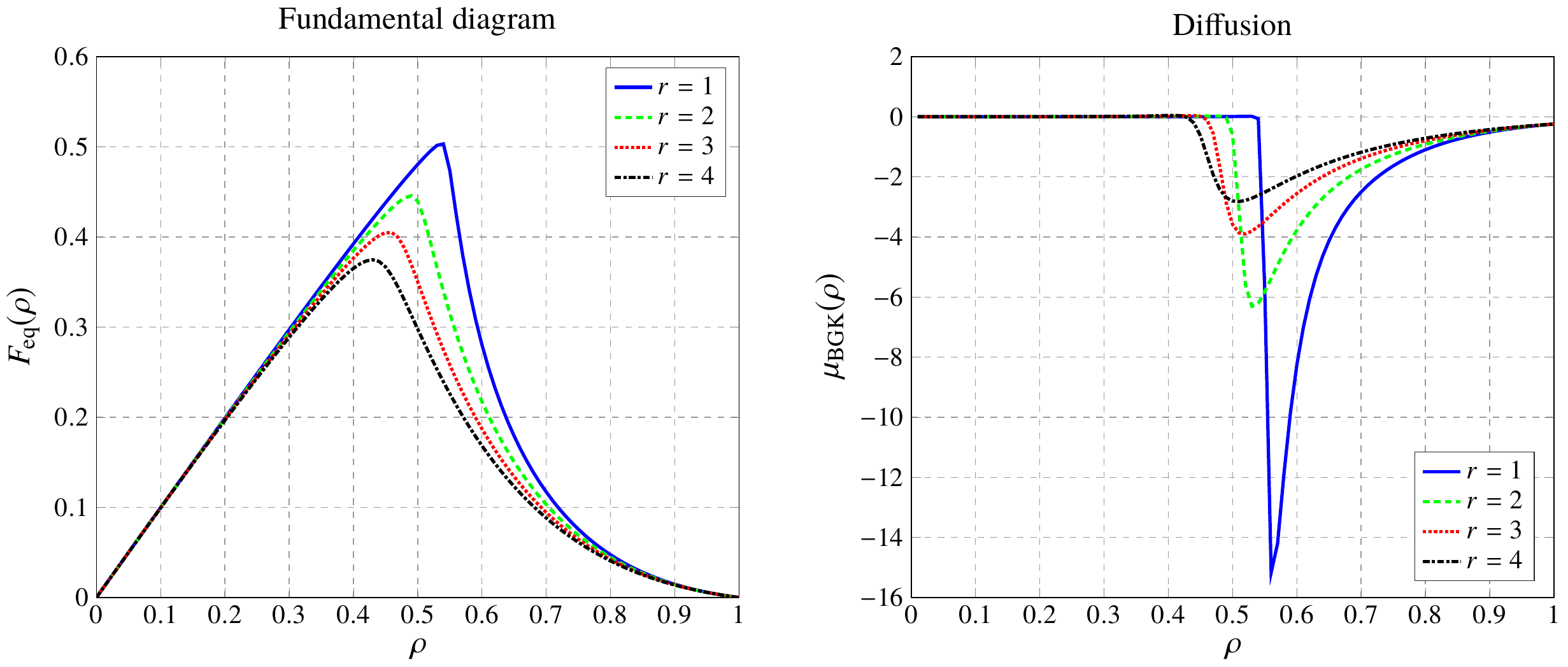}
	\caption{The right panel shows the sign of the diffusion coefficient~\eqref{eq:diffBGK} for the BGK model~\eqref{eq:BGK} with the corresponding equilibrium distribution in the left panel.\label{fig:diffBGK}}
\end{figure}

We  analyze the validity of this result for the model in~\cite{TesiRonc}. In Figure~\ref{fig:diffBGK} we investigate the sign of the diffusion coefficient~\eqref{eq:diffBGK} in the case of the equilibrium distribution. Those distributions are computed numerically for the spatially homogeneous kinetic model~\eqref{eq:generalKineticModel}-\eqref{eq:collisionKernel}. Again, we use $48$ discrete speeds, a fixed value of the acceleration parameter $\Delta_a=\frac{V_M}{4}=\frac14$ and several values of the uncertainty $\Delta_b$ such that $r=\frac{\Delta_a}{\Delta_b}=1,2,3,4$. We observe that $\mu(\rho)_\text{BGK} \geq 0$ in the regime where the flux is increasing, while $\mu(\rho)_\text{BGK} < 0$ in the regime where the flux is decreasing. Increasing the uncertainty on the over-braking, the model becomes ``less'' unstable. In fact, the diffusion becomes larger but still negative. This may serve as explanation of  the growth of perturbations in the density, numerically observed in the top of Figure~\ref{fig:VariableEpsilon}.

In view of the results provided by the Chapman-Enskog analysis, we state the following definition. 

\begin{definition} \label{def:stability}
	A mathematical model for traffic flow is said to be stable if its Chapman-Enskog expansion provides $\mu(\rho)\geq 0$, $\forall\,\rho\in[0,\rho_M]$, weakly-unstable if $\mu(\rho)<0$ on an interval $(\rho_1,\rho_2)$ properly contained in $[0,\rho_M]$ and unstable if $\mu(\rho)<0$ on an interval $(\rho_1,\rho_2)$ in which either $\rho_1=0$ or $\rho_2=\rho_M$.
\end{definition}

The definition of a weakly-unstable model is a consequence of the experimental observation in~\cite{HPRV_BGK2020,SeiboldFlynnKasimovRosales2013} that if $\mu(\rho)<0$ on an interval $(\rho_1,\rho_2)$ properly contained in $[0,\rho_M]$, then the backward propagating waves in dense traffic remain bounded, because, when the oscillations reach $\rho=\rho_1$ and $\rho=\rho_2$, they fall in the diffusive region and they are damped. This  leads to weak instabilities that in turn can be regarded as models for stop-and-go waves.

Concerning the concept of stability of Definition~\ref{def:stability}, the discrete BGK model for traffic introduced in~\cite{BorscheKlar2018} can be either stable or unstable. The model is characterized by non-local terms and with a suitable choice of the headway parameter the diffusion coefficient in the Chapman-Enskog expansion is positive on $[0,\rho_M]$. As observed in~\cite{SeiboldFlynnKasimovRosales2013}, this is not desirable in a model for traffic flow since it would not allow to reproduce non-equilibrium phenomena, such as stop-and-go waves.

\subsection{The Aw-Rascle and Zhang model}

The  Aw-Rascle and Zhang (ARZ) model will be considered in view of the stability analysis following~\cite{HPRV_BGK2020}. We will show that it is weakly-unstable. This justifies the derivation of a new BGK-type model in Section~\ref{sec:4}. The following result was already mentioned and analyzed in~\cite{SeiboldFlynnKasimovRosales2013}. 

The ARZ model reads in primitive variables as
\begin{equation} \label{eq:modelARZnoncons}
\begin{aligned}
\partial_t \rho(x,t) + \partial_x (\rho(x,t) u(x,t)) &= 0\\
\partial_t \big(u(x,t)+h(\rho)\big) + u(x,t) \partial_x \big(u(x,t)+h(\rho)\big) &= \frac{1}{\varepsilon} (U_{\text{eq}}(\rho) - u(x,t)).
\end{aligned}
\end{equation}
where $u$ is the macroscopic speed of the flow and 
the function $h=h(\rho)$ is a strictly increasing function of the density and it is called hesitation function or traffic pressure. The quantity $\varepsilon$ is a time which rules the relaxation speed of the velocity $u$ to the equilibrium speed $U_{\text{eq}}(\rho)$ which is a given function of the density. Here $U_\text{eq}$ is not necessarily  given by ~\eqref{eq:macroQuantitiesEq}.

System~\eqref{eq:modelARZnoncons} can be understood as a relaxation system~\cite{Jin95therelaxation} converging towards the conservation law given by the Lighthill-Whitham~\cite{lighthill1955PRSL} and Richards~\cite{richards1956OR} model in the limit $\varepsilon \to 0$.
If $\varepsilon$ is small, but not vanishing, \eqref{eq:modelARZnoncons}
approaches the advection-diffusion equation~\eqref{eq:advdiff}
where the diffusion coefficient $\mu(\rho)$ is given by
\begin{equation} \label{eq:diffARZeulnoncons}
	\mu_{\text{ARZ}}(\rho) = -\rho(x,t)^2 U^\prime_{\text{eq}}(\rho) \big( U^\prime_{\text{eq}}(\rho) + h^\prime(\rho) \big).
\end{equation}
This result is again obtained via Chapman-Enskog expansion, by considering a first-order expansion of the speed $u=U_{\text{eq}}(\rho) + \varepsilon u_1$ around the equilibrium velocity function $U_{\text{eq}}(\rho)$. The condition $\mu(\rho) > 0$ provides the so-called sub-characteristic condition~\cite{Chen92hyperbolicconservation,Jin95therelaxation}.
For the ARZ model $\mu(\rho)>0$ is satisfied if
\begin{equation} \label{eq:subchCondRed}
0>U_{\text{eq}}^\prime(\rho) > -h^\prime(\rho).
\end{equation}
We  stress the fact that condition~\eqref{eq:subchCondRed} strongly 
restricts the possible choice of $U_\text{eq}$ and $h$, which can be chosen in order to make the model weakly-unstable.

\section{The modified formulation of the BGK approximation in traffic flow} \label{sec:4}

The derivation of the modified BGK-type equation for traffic flow  is shortly summarized and we refer to~\cite{HPRV_BGK2020} for a thorough discussion.
The model is derived via mesoscopic limit of the microscopic follow-the-leader (FTL) and Bando model. We recall that the FTL-Bando model is proved to converge to the ARZ model in the macroscopic limit, both in one-dimension~\cite{aw2002SIAP} and two-dimensions~\cite{HertyMoutariVisconti2018}. Therefore, the second-order system of moments of the new BGK model has also the property of representing a mesoscopic formulation of the class of second-order ARZ-type macroscopic models. As a consequence the feature of an ARZ-type model of having a negative diffusion coefficient in a small density regime is automatically obtained also for the new BGK-type equation.

\subsection{BGK-type model derived from the FTL-Bando model} \label{sec:microFTL}

Let $(x_i,v_i)$ be the microscopic states, position and velocity, of vehicle $i$. The follow-the-leader and Bando model is
\begin{equation} \label{eq:modelFTLw}
\begin{aligned}
\dot{x}_i &= v_i = w_i - p(\rho_i)\\
\dot{w}_i & = \frac{1}{\varepsilon} ( U_\text{eq}(\rho_i) + p(\rho_i) - w_i ).
\end{aligned}
\end{equation}
where $w_i := v_i + p(\rho_i)$, and the function $p = p(\rho_i)$ is the so-called traffic pressure. We assume that $p$ satisfies $p(\rho) \geq 0$, $p^\prime(\rho)>0$ and
$$
\frac{\mathrm{d}}{\mathrm{d}t} p(\rho_i) = - K(x_i,x_{i+1},v_i,v_{i+1})
$$
where $K$ is a term describing the interactions among vehicles. In the classical FTL model
$$
K(x_i,x_{i+1},v_i,v_{i+1}) = C_\gamma \frac{v_{i+1}-v_i}{(x_{i+1}-x_i)^{\gamma+1}},
$$
where the constants $C_\gamma>0$ and $\gamma>0$ are given parameters.
However, we consider the case of a general function $K$. The introduction of the quantity $w_i$ allows us to rewrite the classical Bando model as a relaxation step~\eqref{eq:modelFTLw}.

Let now $g = g(x,w,t) : \mathbb{R} \times W \times \mathbb{R}^+ \to \mathbb{R}^+$ be the kinetic distribution function with respect to the desired speed $w$, which is assumed to be the speed that drivers want to keep in ``optimal'' situations.
We define $W:=[w_{\min},+\infty)$ the space of the microscopic desired speeds where $w_{\min}>0$ may be interpreted as the minimum speed limit in free-flow conditions. The macroscopic density, i.e. the number of vehicles per unit length, at time $t$ and position $x$ is defined by
\begin{equation} \label{eq:densityg}
\rho(x,t) := \int_W g(x,w,t) \dd w,
\end{equation}
and we define the macroscopic quantity
\begin{equation} \label{eq:defq}
q(x,t) := \int_W w g(x,w,t) \dd w.
\end{equation}

The derivation of the evolution equation for the kinetic distribution $g=g(x,w,t)$ is performed by reformulating the microscopic particle model~\eqref{eq:modelFTLw} in a probabilistic interpretation and  allowing a relaxation  towards a desired distribution $M_g = M_g(w;\rho)$, as in~\cite[Section 4.2.2]{PareschiToscaniBOOK}. The distribution $M_g$ has to fulfill the requirement
$$
\int_W M_g(w;\rho) \dd w = \rho(x,t),
$$
and additionally
\begin{equation} \label{eq:MgRelation}
\frac{1}{\rho(x,t)} \int_W w M_g(w;\rho) \dd w = U_\text{eq}(\rho) + p(\rho).
\end{equation}

According to~\cite[Section 4.2.2]{PareschiToscaniBOOK} and~\cite{HPRV_BGK2020}, it is possible to show that $g$ solves
\begin{equation} \label{eq:kineticW}
\partial_t g(x,w,t) + \partial_x \big[ (w-p(\rho)) g(x,w,t) \big] = \frac{1}{\varepsilon} \left( M_g(w;\rho) - g(x,w,t) \right)
\end{equation}
This equation  is still a BGK-type equation since the collision kernel is linear and describes the relaxation of $g$ towards a given distribution $M_g$ parameterized by the density $\rho$. For a detailed derivation of~\eqref{eq:kineticW} we refer to~\cite{HPRV_BGK2020}. 

It is important to point-out that, compared to classic kinetic theory, this approach is different in the sense that $M_g$ is an ``equilibrium distribution'' with a modified  microscopic velocity. Thanks to~\eqref{eq:MgRelation}, $M_g$ is imposed a-priori but it is still based on the knowledge of the classical Maxwellian $M_f$, which is related to the classical concept of microscopic velocity, by means of $U_\text{eq}(\rho) := \frac{1}{\rho} \int v M_f \dd v$. In other words, $M_f$ is not imposed a-priori (and so $U_\text{eq}(\rho)$ and consequently $M_g$), but the equilibrium distribution $M_f$ is the one obtained by the modeling of microscopic interactions of the spatially homogeneous kinetic model.  Any Maxwellian $M_f$ of a kinetic models for traffic can be used to define $M_g$ and the BGK model~\eqref{eq:kineticW}. Here,  we study the Maxwellian $M_f$ provided by~\cite{TesiRonc}.

\subsection{Chapman-Enskog expansion of the modified BGK-model}

We perform a Chapman-Enskog expansion for the model~\eqref{eq:kineticW}. We consider a first-order perturbation of $g$ as
$$
g(x,w,t) = M_g(w;\rho) + \varepsilon g_1(x,w,t), \quad \text{with} \ \int_W g_1(x,w,t) \dd w = 0
$$
and define $F_\text{eq}(\rho) = \rho U_\text{eq}(\rho)$.
Then, it is possible to show, cf.~\cite{HPRV_BGK2020}, that the BGK-type equation~\eqref{eq:kineticW} solves the advection-diffusion equation~\eqref{eq:advdiff}
with
\begin{equation} \label{eq:kineticDiffCoeff}
\mu(\rho) = -F^\prime_\text{eq}(\rho)^2 + \int_V v^2 \partial_\rho M_f(v;\rho) \dd v - \rho p^\prime(\rho) F^\prime_\text{eq}(\rho) + F_\text{eq}(\rho) p^\prime(\rho).
\end{equation}
Observe that, compared to~\eqref{eq:diffBGK}, the diffusion coefficient~\eqref{eq:kineticDiffCoeff} contains two additional terms which depend on the function $p(\rho)$. Therefore, it is possible, for a given distribution $M_f$, to find a suitable $p(\rho)$ such that $\mu(\rho) > 0$ also in the congested regime. In particular, it is possible to find $p(\rho)$ in order to guarantee that the model is weakly-unstable. Recall that $\mu_\text{BGK}(\rho)$ given in~\eqref{eq:diffBGK} was unconditionally negative in the congested phase of traffic for the classical BGK model~\eqref{eq:BGK}.

\begin{figure}[t!]
	\centering
	\includegraphics[width=\textwidth]{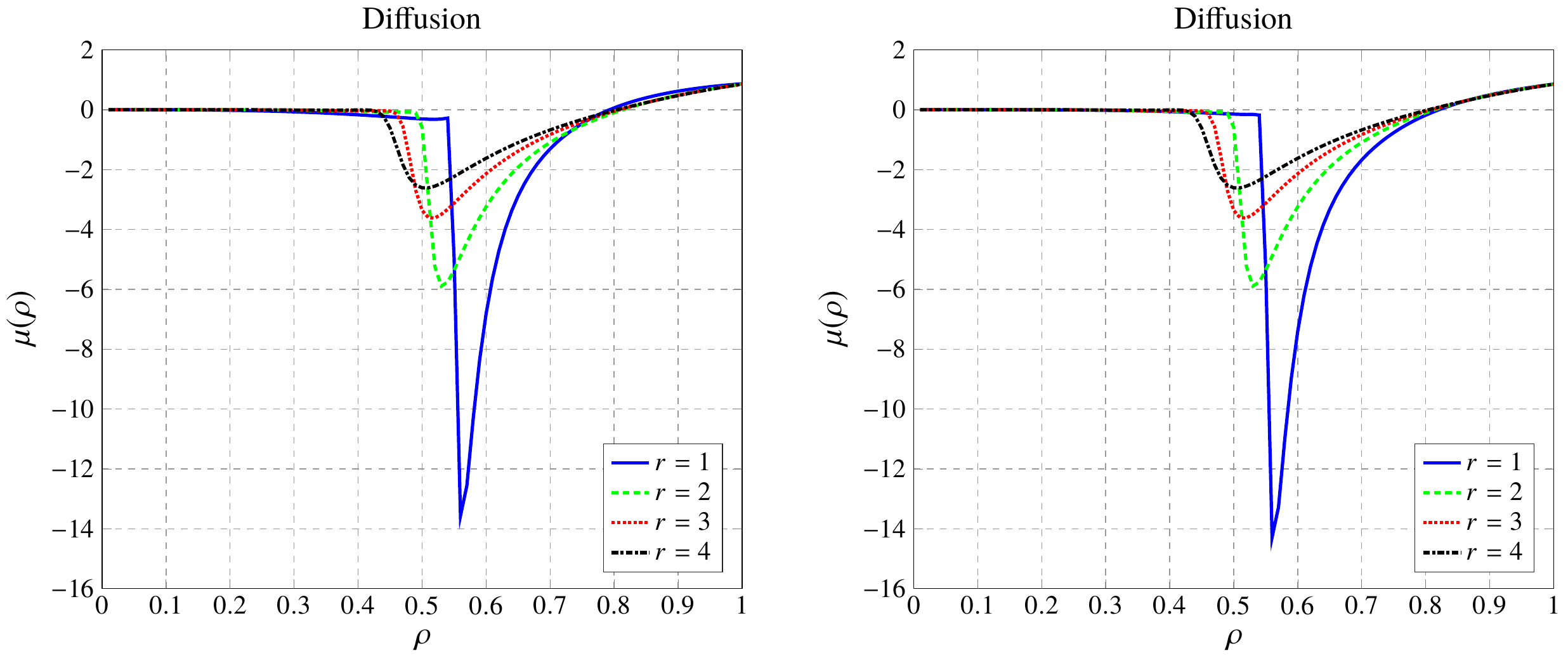}
	\caption{Diffusion coefficient~\eqref{eq:newDiffusion} for $p(\rho)=\frac32\rho^2$ (left) and $p(\rho)=\rho^3$ (right).}
			\label{fig:diffusionNewBGK}
\end{figure}

Setting $F_\text{eq}(\rho) = \rho U_\text{eq}(\rho)$ the second two terms of the diffusion coefficient~\eqref{eq:kineticDiffCoeff} can be written in terms of the equilibrium speed function as
\begin{equation} \label{eq:newDiffusion}
\mu(\rho) = -F^\prime_\text{eq}(\rho)^2 + \int_V v^2 \partial_\rho M_f(v;\rho) \mathrm{d}v - \rho^2 p^\prime(\rho) U_\text{eq}^\prime(\rho).
\end{equation}
Therefore, $\mu(\rho)=\mu_\text{BGK}(\rho)+C(\rho)$ where $C(\rho) = - \rho^2 p^\prime(\rho) U_\text{eq}^\prime(\rho)\geq 0$ since $p$ and $U_\text{eq}$ are an increasing and a non-increasing function of the density, respectively. This means that, for  $C(\rho)$ sufficiently large, the additional term  yields a negative diffusion coefficient~\eqref{eq:diffBGK}  in a bounded domain  contained in $[0,\rho_M]$. In Figure~\ref{fig:diffusionNewBGK} we numerically show this result for the case of the homogeneous kinetic model in~\cite{TesiRonc}. We consider the Maxwellian computed numerically with $48$ discrete speeds, $\Delta_a=\frac14$ and $\Delta_b=\frac{\Delta_a}{r}$, $r=1,2,3,4$. The pressure function is chosen as $p(\rho) = \frac32 \rho^2$ (left panel) and $p(\rho) = \rho^3$ (right panel).

\section{Conclusions and future perspectives} \label{sec:conslusion}

In this work we have focused on the formulation of kinetic models for vehicular traffic flow which reproduce backward propagating waves in dense traffic. The underlying kinetic model is the one introduced in~\cite{TesiRonc}. Backward traveling waves have been observed by defining an interaction rate that is a function of the density and its derivative. 

A stability analysis of the waves in dense traffic regimes has been performed on the BGK-type approximation, in the limit of small interaction rates. We have shown that the model leads to an advection-diffusion equation with a negative diffusion coefficient in the whole congested regime, therefore producing an unbounded growth of dense waves in time. This  justified to reconsider the results of \cite{TesiRonc} in the framework of a novel BGK formulation~\cite{HPRV_BGK2020}. Finally, the formulation allows to have a weakly-unstable model with results that show the existence of stop-and-go waves.

\section*{Acknowledgments}

The research of M. Herty and G. Visconti is funded by the Deutsche Forschungsgemeinschaft (DFG, German Research Foundation) under Germany's Excellence Strategy -- EXC-2023 Internet of Production -- 390621612 as well as by DFG HE5386/13.  

G. Puppo and G. Visconti acknowledge also support from GNCS (Gruppo Nazionale per il Calcolo Scientifico) of INdAM (Istituto Nazionale di Alta Matematica), Italy.

\bibliographystyle{plain}
\bibliography{AppuntiTraBib}

\begin{thebibliography}{10}

\bibitem{aw2002SIAP}
A.~Aw, A.~Klar, T.~Materne, and M.~Rascle.
\newblock Derivation of continuum traffic flow models from microscopic
  follow-the-leader models.
\newblock {\em SIAM J. Appl. Math.}, 63(1):259--278, 2002.

\bibitem{aw2000SIAP}
A.~Aw and M.~Rascle.
\newblock Resurrection of ``second order'' models of traffic flow.
\newblock {\em SIAM J. Appl. Math.}, 60(3):916--938 (electronic), 2000.

\bibitem{Bando1995}
M.~Bando, K.~Hasebe, A.~Nakayama, A.~Shibata, and Y.~Sugiyama.
\newblock Dynamical model of traffic congestion and numerical simulation.
\newblock {\em Phys. Rev. E}, 51(2):1035--1042, 1995.

\bibitem{BGK1954}
P.~L. Bhatnagar, E.~P. Gross, and M.~Krook.
\newblock A {M}odel for {C}ollision {P}rocesses in {G}ases. {I}. {S}mall
  {A}mplitude {P}rocesses in {C}harged and {N}eutral {O}ne-{C}omponent
  {S}ystems.
\newblock {\em Phys. Rev.}, 94(3):511--525, 1954.

\bibitem{BorscheKlar2018}
R.~Borsche and A.~Klar.
\newblock A nonlinear discrete velocity relaxation model for traffic flow.
\newblock {\em SIAM J. Appl. Math.}, 78(5):2891--2917, 2018.

\bibitem{Chen92hyperbolicconservation}
G.-q. Chen, C.~D. Levermore, and T.-P. Liu.
\newblock Hyperbolic conservation laws with stiff relaxation terms and entropy.
\newblock {\em Comm. Pure Appl. Math}, 47:787--830, 1992.

\bibitem{DimarcoPareschi2013}
G.~Dimarco and L.~Pareschi.
\newblock {Asymptotic Preserving Implicit-Explicit Runge-Kutta methods for non
  linear kinetic equations}.
\newblock {\em SIAM J. Num. Anal}, 51:1064--1087, 2013.

\bibitem{FermoTosin13}
L.~Fermo and A.~Tosin.
\newblock A fully-discrete-state kinetic theory approach to modeling vehicular
  traffic.
\newblock {\em SIAM J. Appl. Math.}, 73(4):1533--1556, 2013.

\bibitem{FTL1961}
D.~Gazis, R.~Herman, and R.~Rothery.
\newblock Nonlinear follow-the-leader models of traffic flow.
\newblock {\em Oper. Res.}, 9(4):545--567, 1961.

\bibitem{HertyMoutariVisconti2018}
M.~Herty, S.~Moutari, and G.~Visconti.
\newblock Macroscopic modeling of multilane motorways using a two-dimensional
  second-order model of traffic flow.
\newblock {\em SIAM J. Appl. Math.}, 78(4):2252--2278, 2018.

\bibitem{HPRV_BGK2020}
M.~Herty, G.~Puppo, S.~Roncoroni, and G.~Visconti.
\newblock The {BGK} approximation of kinetic models for traffic.
\newblock {\em Kinet. Relat. Models}, 2020.
\newblock In press.

\bibitem{Jin95therelaxation}
S.~Jin and Z.~Xin.
\newblock The relaxation schemes for systems of conservation laws in arbitrary
  space dimensions.
\newblock {\em Comm. Pure Appl. Math}, 48:235--277, 1995.

\bibitem{klar1997Enskog}
A.~Klar and R.~Wegener.
\newblock Enskog-like kinetic models for vehicular traffic.
\newblock {\em J. Stat. Phys.}, 87:91, 1997.

\bibitem{lighthill1955PRSL}
M.~J. Lighthill and G.~B. Whitham.
\newblock On kinematic waves. {II}. {A} theory of traffic flow on long crowded
  roads.
\newblock {\em Proc. Roy. Soc. London. Ser. A.}, 229:317--345, 1955.

\bibitem{PareschiToscaniBOOK}
L.~Pareschi and G.~Toscani.
\newblock {\em Interacting {M}ultiagent {S}ystems. {K}inetic equations and
  {M}onte {C}arlo methods}.
\newblock Oxford University Press, 2013.

\bibitem{PgSmTaVg3}
G.~Puppo, M.~Semplice, A.~Tosin, and G.~Visconti.
\newblock Analysis of a multi-population kinetic model for traffic flow.
\newblock {\em Commun. Math. Sci.}, 15(2):379--412, 2017.

\bibitem{PSTV2}
G.~Puppo, M.~Semplice, A.~Tosin, and G.~Visconti.
\newblock Kinetic models for traffic flow resulting in a reduced space of
  microscopic velocities.
\newblock {\em Kinet. Relat. Mod.}, 10(3):823--854, 2017.

\bibitem{richards1956OR}
P.~I. Richards.
\newblock Shock waves on the highway.
\newblock {\em Operations Res.}, 4:42--51, 1956.

\bibitem{TesiRonc}
Sebastiano Roncoroni.
\newblock Kinetic modelling of vehicular traffic flow.
\newblock Technical report, Universit\`a degli Studi dell'Insubria, 2017.
\newblock Master Thesis.

\bibitem{SeiboldFlynnKasimovRosales2013}
B.~Seibold, M.~R. Flynn, A.~R. Kasimov, and R.~R. Rosales.
\newblock Constructing set-valued fundamental diagrams from jamiton solutions
  in second order traffic models.
\newblock {\em Netw. Heterog. Media}, 8(3):745--772, 2013.

\bibitem{Zhang2002}
H.~M. Zhang.
\newblock A non-equilibrium traffic model devoid of gas-like behavior.
\newblock {\em Transport. Res. B-Meth.}, 36(3):275--290, 2002.

\end{thebibliography}

\end{document}